\documentclass[12pt]{article} \newcommand{\ifcmp}[1]{} \newcommand{\ifarxiv}[1]{#1}

\usepackage{amsmath,amssymb,amsfonts,graphicx}
\ifarxiv{\usepackage{amsthm}}
\usepackage[pdftitle={Conformal radii for conformal loop ensembles},
            pdfauthor={Oded Schramm, Scott Sheffield, and David B. Wilson},
            bookmarks=true,bookmarksopen=true,bookmarksopenlevel=2,unicode]{hyperref}

\ifcmp{
\journalname{Communications in Mathematical Physics}
\date{Received: 20 September 2007 / Accepted: 10 November 2008}
\communicated{Michael Aizenman}
}
\ifarxiv{
\date{}
}

\newcommand{\sref}[1]{\hyperref[#1]{Section~\ref*{#1}}}
\newcommand{\aref}[1]{\hyperref[#1]{Appendix~\ref*{#1}}}
\newcommand{\lref}[1]{\hyperref[#1]{Lemma~\ref*{#1}}}
\newcommand{\tref}[1]{\hyperref[#1]{Theorem~\ref*{#1}}}
\newcommand{\cref}[1]{\hyperref[#1]{Corollary~\ref*{#1}}}
\newcommand{\fref}[1]{\hyperref[#1]{Figure~\ref*{#1}}}
\newcommand{\pref}[1]{\hyperref[#1]{Proposition~\ref*{#1}}}
\newcommand{\eref}[1]{\hyperref[#1]{Equation~\ref*{#1}}}

\newcommand{\dref}[1]{\hyperref[#1]{Diffusion~\ref*{#1}}}

\def\clap#1{\hbox to 0pt{\hss#1\hss}}

 \newcommand{\MRhref}[2]{\href{http://www.ams.org/mathscinet-getitem?mr=#1}{MR#2}}

\makeatletter
\def\@strippedMR{}
\def\@scanforMR#1#2#3\endscan{%
  \ifx#1M\ifx#2R\def\@strippedMR{#3}%
  \else\def\@strippedMR{#1#2#3}%
  \fi\fi}
\def\@rst #1 #2other{#1}
\newcommand\MR[1]{\relax\ifhmode\unskip\spacefactor3000 \space\fi
  \@scanforMR#1\endscan
  \MRhref{\expandafter\@rst \@strippedMR other}{\@strippedMR}}
\newcommand\MRs[1]{\relax\ifhmode\unskip\spacefactor3000 \space\fi
  \@scanforMR#1\endscan
  \MRhref{\@strippedMR}{\@strippedMR}}
\makeatother

\ifarxiv{
 \setlength{\textwidth}{6.5 in}
 \setlength{\textheight}{8.25in}
 \setlength{\oddsidemargin}{0in}
 \setlength{\topmargin}{0in}
 \addtolength{\textheight}{.8in}
 \addtolength{\voffset}{-.5in}
\newtheorem{theorem}{Theorem}

\newtheorem{lemma}{Lemma}

\newtheorem{proposition}{Proposition}
}
\newcommand{\Z}{\mathbb{Z}}
\newcommand{\R}{\mathbb{R}}
\renewcommand{\C}{\mathbb{C}}

\newcommand{\N}{\mathbb{N}}
\newcommand{\E}{\mathbb{E}}
\newcommand{\D}{\mathbb{D}}
\newcommand{\eps}{\varepsilon}

\newcommand{\dist}{\operatorname{dist}}
\def\SLEkk#1/{$\mathrm{SLE}_{#1}$}
\def\SLEk/{\SLEkk{\kappa}/}
\def\SLEtwo/{\SLEkk2/}
\def\SLE/{$\mathrm{SLE}$}
\def\CLEkk#1/{$\mathrm{CLE}_{#1}$}
\def\CLEk/{\CLEkk{\kappa}/}
\def\CLEtwo/{\CLEkk2/}
\def\CLE/{$\mathrm{CLE}$}
\def\Ito{It\^o}
\newcommand{\crad}{\operatorname{CR}}

\renewcommand{\Re}{\operatorname{Re}}

\newcommand{\RY}{MR2000h:60050}
\newcommand{\BE}{bateman-erdelyi:vol1}
\ifarxiv{\pagestyle{myheadings} \markright{Conformal radii for conformal loop ensembles \hfil Schramm, Sheffield, \& Wilson}}
\begin{document}
\title{Conformal radii for conformal loop ensembles}
\ifcmp{
\author{Oded Schramm\inst{1} \and
Scott Sheffield\inst{2}\thanks{Partially supported by NSF grant DMS0403182.} \and
 David B. Wilson\inst{1}}}
\ifarxiv{
\author{Oded Schramm\thanks{Microsoft Research.} \and
Scott Sheffield\thanks{Courant Institute and Massachusetts Institute of Technology.  Partially supported by NSF grant DMS0403182.} \and
 David B. Wilson\footnotemark{1}}}
\ifcmp{\institute{Microsoft Research \and Courant Institute and Massachusetts Institute of Technology}}
\maketitle

\begin{abstract}
The conformal loop ensembles \CLEk/, defined for $8/3 \leq \kappa
\leq 8$, are random collections of loops in a planar domain which
are conjectured scaling limits of the $O(n)$ loop models. We
calculate the distribution of the conformal radii of the nested
loops surrounding a deterministic point.  Our results agree with
predictions made by Cardy and Ziff and by Kenyon and Wilson for the
$O(n)$ model.  We also compute the expectation dimension of the
\CLEk/ \textbf{gasket}, which consists of points not surrounded by any
loop, to be
$$ 2-\frac{(8-\kappa) (3\kappa-8)}{32\kappa}\,,$$ which agrees with the fractal dimension
given by Duplantier for the $O(n)$ model gasket.
\end{abstract}

\section{Introduction}

The conformal loop ensembles \CLEk/, defined for all $8/3 \leq
\kappa \leq 8$, are random collections of loops in a simply
connected planar domain~$D\subsetneqq \C$.  They were defined and constructed from
branching variants of \SLEk/ in \cite{Sheffield}, where they were
conjectured to be the scaling limits of various random loop models
from statistical physics, including the so-called $O(n)$ loop models
with
\begin{equation} \label{k-n}
  n=-2\cos(4\pi/\kappa)\,,
\end{equation}
see e.g.\ \cite{MR2065722} for an exposition.  This paper is a
sequel to \cite{Sheffield}.  We will state the results about \CLEk/
from \cite{Sheffield} that we need for this paper (namely
Propositions~\ref{p.recursive} and \ref{p.bklaw}), but we will not
repeat the definition of \CLEk/ here.

When $8/3 < \kappa < 8$, \CLEk/ is almost surely a countably
infinite collection of loops. \CLEkk 8/ is a single space-filling
loop almost surely and \CLEkk {8/3}/ is almost surely empty. \CLEkk
6/ is the scaling limit of the cluster boundaries of critical site
percolation on the triangular lattice \cite{camia-newman:full-limit}
\cite{MR1851632} \cite{camia-newman:proof-converge}.
We will henceforth assume $8/3 < \kappa < 8$.

For each $z \in D$, we inductively
define $L_k^z$ to be the outermost loop surrounding~$z$ when the loops
$L_1^z, \ldots, L_{k-1}^z$ are removed (provided such a loop
exists). For each deterministic~$z \in D$, the loops~$L_k^z$ exist
for all~$k \geq 1$ with probability one.  Define~$A_0^z = D$ and let
$A_k^z$ be the component of $D \setminus L_k^z$ that contains~$z$.
The \textbf{conformal gasket} is the random closed set~$\Gamma$
consisting of points that are not surrounded by any loop of an
instance of \CLEk/, i.e.\ the set of points for which~$L_1^z$ does
not exist.

If $D$ is a simply connected planar domain and~$z \in D$, the {\bf
conformal radius of~$D$ viewed from~$z$} is defined to be
$\crad(D,z) := |g'(z)|^{-1}$, where $g$ is any conformal map from~$D$
to the unit disk~$\D$ that sends $z$ to~$0$.  The following is
immediate from the construction in \cite{Sheffield}:

\begin{proposition} \label{p.recursive}
Let $D$ be a simply connected bounded planar domain,
and consider a \CLEk/ on~$D$ for some $8/3 < \kappa < 8$.  Then
$\Gamma$ is almost surely the closure of the set of points that lie
on an outermost loop (i.e., a loop of the form~$L_1^z$ for some~$z$).
Conditioned on the outermost loops, the law of the remaining
loops is given by an independent \CLEk/ in each component of $D
\setminus \Gamma$.  For $z\in D$ and $k=1,2,3,\dots$, define
$$B_k^z := \log \crad(A_{k-1}^z, z) - \log \crad(A_k^z,z)\,.$$
For any fixed~$z$, the $B_k^z$'s are i.i.d.\ random variables.
\end{proposition}

Various authors in the physics literature have used heuristic
arguments (based on the so-called Coulomb gas method) to calculate
properties of the scaling limits of statistical physical loop
models, including the $O(n)$ models, based on certain conformal
invariance hypotheses of these limits. Although the scaling limits
of the $O(n)$ models have not been shown to exist, there is strong
evidence that if they do exist they must be \CLEk/.  (For example,
there is heuristic evidence that any scaling limit of the $O(n)$
models should be in some sense conformally invariant; it is shown in
\cite{Sheffield} \cite{SheffieldWerner:Markov} \cite{SheffieldWerner:Soup}
 that any random loop
ensemble satisfying certain hypotheses including conformal
invariance and a Markov-type property must be a \CLEk/.)  It is
therefore natural to interpret these calculations as predictions
about the behavior of the \CLEk/.

Cardy and Ziff \cite{CZ} predicted and experimentally verified the
expected number of loops surrounding a point in the $O(n)$ model,
which, in light of \eqref{k-n}, may be interpreted as a prediction of
the expectation of~$B^k_z$:
\begin{equation}
  \label{expected-radii}
 \frac{1}{ \E[B_k^z]} = \frac{(\kappa/4-1)\cot(\pi(1-4/\kappa))}{\pi} \,.
\end{equation}
Kenyon and Wilson \cite{KW} went further and predicted the
distribution of~$B_k^z$, giving its moment generating function
\begin{equation} \label{mgf}
  \E[\exp(\lambda B_k^z)] = \frac{-\cos(4\pi/\kappa)}
                 {\cos\left(\pi\sqrt{(1-4/\kappa)^2+ 8 \lambda /\kappa}\right)}
\,,
\end{equation}
for $\lambda$ satisfying
$\Re \lambda < 1-\frac 2{\kappa}-\frac {3 \kappa}{32}$,
and density function
\begin{equation} \label{pdf}
\frac{d}{dx} \Pr[B_k^z<x] = \frac{-\kappa \cos(4\pi/\kappa)}{4\pi}
 \sum_{j=0}^\infty (-1)^j (j+1/2)
                  \exp\left[-\frac{(j+1/2)^2-(1-4/\kappa)^2}{8/\kappa}x\right]\,.
\end{equation}
The main result of this paper is \tref{Bklaw}, which confirms
these predictions.  In the special case~$\kappa=6$, this prediction
for the law of~$B_k^z$ was independently confirmed by Dub\'edat
\cite{dubedat}.

\begin{theorem} \label{Bklaw}
  Let $f_\kappa$ denote the density function for the first time that a
  standard Brownian motion started at~$0$ exits the interval $(-
  2\pi/\sqrt{\kappa}, 2\pi/\sqrt{\kappa})$.  Then for $8/3<\kappa<8$,
  the density function for~$B_k^z$ is
\begin{equation} \label{e.density2}
   \frac{d}{dx} \Pr[B_k^z<x] =
  -f_\kappa(x) \cos(4\pi/\kappa) \exp\left[\frac{(\kappa-4)^2}{8 \kappa} x\right]\,.
\end{equation}
\end{theorem}
The equivalence of the formulations \eqref{pdf} and
\eqref{e.density2}, and the fact that they imply \eqref{mgf}, follows
from a calculation of Ciesielski and Taylor, who showed that the
exit-time distribution of a Brownian motion from the center of a
1-dimensional ball of radius~$r$ has a moment generating function of
$1/\cos\sqrt{2 r^2 \lambda}$, and who gave two series expansions (one
in powers of~$e^{-x}$ and the other in powers of~$e^{-1/x}$) for its
density function \cite[Theorem~2 and eq.~2.22]{ciesielski-taylor} (see
also \cite[eqs.~1.3.0.1 and 1.3.0.2]{borodin-salminen}).  Since the
Fourier transform is invertible on~$L^2(\R)$, the equivalence of
\eqref{mgf} and \eqref{pdf} follows by considering the moment
generating function restricted to the imaginary line~$\Re\lambda=0$.

Duplantier \cite{D:vesicle} predicted the fractal dimension of the
gasket~$\Gamma$ associated with the $O(n)$ model to be
$$
     \frac{3\kappa}{32} + 1 + \frac{2}{\kappa}\,,
$$
where as usual $n=-2\cos(4\pi/\kappa)$.  We partially confirm this
prediction by giving the expectation dimension of the gasket
associated with \CLEk/.
The \textbf{expectation dimension} of a random bounded set~$A$ is defined to be
$$
D_{\E}(A) = \lim_{\eps\to 0} \frac{\log\E[\text{minimal number of balls of radius~$\eps$ required to cover~$A$}]}{|\log\eps|}\,,
$$
provided the limit exists.  The expectation dimension upper bounds
the Hausdorff dimension.
\begin{theorem} \label{t.dimension}
Let $\Gamma$ be the gasket  of a \CLEk/ in the unit disk
with $\kappa\in(8/3,8)$.  Then
$$
\E[\text{minimal number of balls of radius~$\eps$ required to cover~$\Gamma$}] \asymp \left(\frac1\eps\right)^{\frac{3\kappa}{32} + 1 + \frac{2}{\kappa}}\,.
$$
In particular, the expectation dimension of the gasket~$\Gamma$ is
$
     \frac{3\kappa}{32} + 1 + \frac{2}{\kappa}.
$
\end{theorem}
Here, $\asymp$ denotes equivalence up to multiplicative constants.
Lawler, Schramm, and Werner \cite{LSW:one-arm} studied the
percolation gasket (associated with \CLEkk6/), effectively
proving \tref{t.dimension} in the case~$\kappa=6$.
More generally, they studied how long it takes for a radial
\SLEk/ to surround the origin when~$\kappa>4$, and
their results implicitly imply Theorem~2 when $\kappa>4$;
see the remark in \sref{thetat} for further discussion.
\medskip

We conclude our introduction by noting that the gasket dimension
described above plays an important role in the physics literature,
where it is related (at least heuristically) to the exponents of
magnetization and multipoint correlation functions in critical lattice
models.  We briefly describe this connection in the case of the
$q$-state Potts model on the square lattice.  More details and
references are found in \cite{Sheffield}, \cite{Cardy:ADE}, \cite{MR2243761}.

A sample from the $q$-state Potts model on a connected planar graph~$G$
is a random function $\sigma: V \to \{1,2, \ldots, q\}$, where $V$
is the set of vertices of~$G$ and the image values $1,2,\ldots, q$ are
often called spins.  If the boundary vertices
(those on the boundary of the unbounded face) of~$G$ are all assigned
a particular value (say~$b$), then using the standard FK random cluster
decomposition \cite{MR0359655}, one may construct a sample from the Potts model
as follows:

\begin{enumerate}
\item Sample a random subgraph~$G'$ of~$G$ containing all boundary
  edges (edges on the boundary of the unbounded face), with
  probability proportional to $$q^{\text{\# components of $G'$}}
  \left(\frac{p}{1-p}\right)^{\text{\# edges of $G'$}}\,,$$
  where $0<p<1$ is a parameter.
The law of~$G'$ is called the FK random cluster model with parameters~$p$ and~$q$.
  Call the component of~$G'$ which contains the boundary vertices of~$G$ the
\textbf{FK gasket}.

\item Set $\sigma(v) = b$ for each $v$ in the FK gasket, and
  independently assign one of the $q$ states uniformly at random to
  each of the remaining connected components of~$G'$ (assigning all vertices in
  that component the corresponding state).
\end{enumerate}

The ``magnetization'' at an interior vertex~$v$ of~$G$ (i.e., the
probability that $\sigma(v) = j$ minus~$1/q$) is proportional to the
probability that $v$ is in the FK gasket. 
Given distinct vertices~$v$ and~$w$, the covariance of
$\sigma(v)$ and~$\sigma(w)$ is proportional to the probability that
both~$v$ and~$w$ lie in the same component of~$G'$.

We now restrict to the case in which $G$ is a finite piece of the
square grid in the plane and the parameter~$p$ satisfies $p/(1-p)=\sqrt{q}$.
(With this choice of~$p$, the FK model is self-dual and believed to be critical, see e.g., \cite[Chapter~6]{MR2243761}.)

It is shown in \cite{Sheffield} that if $q\leq 4$ and certain
other hypotheses including conformal invariance hold, then the scaling
limit of the set of boundaries between clusters and dual clusters in
the critical FK random cluster models discussed above must be given by \CLEk/ for
the~$\kappa$ satisfying $q=4\cos^2(4\pi/\kappa)$ and $4\leq\kappa\leq 8$.
Assuming these hypotheses, the scaling limit of the discrete gasket is
precisely the continuum \CLEk/ gasket.

A heuristic ansatz is that the law of the critical FK gasket should
have similar properties as the law of the set of squares in a fine
grid which intersect the continuum gasket.  If this heuristic holds,
then when $G$ is a bounded domain intersected with a square grid with
spacing~$\eps$, the magnetization at a vertex~$v$ of macroscopic
distance from the boundary in the discrete model will be on the order
of~$\eps^{2-d}$, where $d$ is the limiting expectation dimension of
the continuum gasket.  Similarly, the covariance between $\sigma(v)$
and~$\sigma(w)$, for two macroscopically separated vertices~$v$ and~$w$,
should be on the order of~$\eps^{2(2-d)}$ (since in the continuum
model, the set of pairs $v$ and~$w$ which lie in the same continuum
spin cluster has dimension~$2d$; see \cite{Sheffield}).

\section{Diffusions and martingales}
\subsection{Reduction to a diffusion problem} \label{thetat}

Let $B_t:[0,\infty)\to \R$ be a standard Brownian motion and let
$\theta_t:[0,\infty) \to [0, 2\pi]$ be a random continuous process
on the interval~$[0,2\pi]$ that is instantaneously reflecting at its
endpoints (i.e., the set $\bigl\{t:\theta_t \in\{0,2\pi\} \bigr\}$ has
Lebesgue measure zero almost surely) and evolves according to the
SDE
\begin{equation}\label{e.mainSDE}d\theta_t = \frac{\kappa-4}{2} \cot(\theta_t/2)\,dt +
\sqrt{\kappa}\, dB_t
\end{equation}
on each interval of time for which $\theta_t \notin \{0,
2\pi\}$. In other words, $\theta_t$ is a random continuous process
adapted to the filtration of~$B_t$ which almost surely satisfies
$$\frac{\partial}{\partial t} (\theta_t - \sqrt{\kappa}\, B_t) =\frac{\kappa-4}{2}
\cot(\theta_t/2)$$ for all~$t$ for which the right hand side is well
defined.  The law of this process is uniquely determined by~$\theta_0$
\cite{Sheffield}, and we also have the following from
\cite{Sheffield}:

\begin{proposition} \label{p.bklaw}
  When $8/3<\kappa<8$, the law of~$B_k^z$ is the same as the law of
  $\inf \{t:\theta_t = 2\pi \}$ for the diffusion~\eqref{e.mainSDE}
  started at~$\theta_0=0$.
\end{proposition}

It is convenient to lift the process~$\theta_t$ so that, rather than
taking values in~$[0,2\pi]$, it takes values in all of~$\R$.  Let
$R:\R \to [0, 2\pi]$ be the piecewise affine map for which $R(x) =
|x|$ when $x \in [-2\pi, 2\pi]$ and $R(4\pi + x) = R(x)$ for all~$x$.
Given~$\theta_t$, we can generate a continuous process~$\tilde
\theta_t$ with $R(\tilde \theta_t) = \theta_t$ in such a way that
for each component~$(t_1, t_2)$ of the set $\{t: \theta_t \notin
2\pi \Z \}$, we independently toss a fair coin to decide whether
$\tilde \theta_t > \tilde \theta_{t_1}$ or $\tilde \theta_t < \tilde \theta_{t_1}$
on that component.  The~$\theta_t$ together with these coin tosses
(for each interval of $\{t: \theta_t \notin 2\pi \Z \}$)
determine~$\tilde \theta_t$ uniquely.

This~$\tilde \theta_t$ is still a solution to \eqref{e.mainSDE}
provided we modify~$B_t$ in such a way that $dB_t$ is replaced with~$-dB_t$
on those intervals for which $d \tilde \theta_t = -
d\theta_t$.  (This modification does not change the law of~$B_t$.)
In the remainder of the text, we will drop the $\tilde \theta_t$
notation and write $\theta_t$ for the lifted process on~$\R$.

\medskip
\noindent
\textbf{Remark: }
A very similar diffusion process was studied by Lawler, Schramm, and Werner \cite{LSW:one-arm}, namely
\begin{equation}\label{e.one-arm}d\theta_t = \cot(\theta_t/2)\,dt +
\sqrt{\kappa}\, dB_t \,,
\end{equation}
which is the same as \eqref{e.mainSDE} but without the factor of
$(\kappa-4)/2$, and they too studied the time for the diffusion to
reach~$\theta_t=2\pi$ when started at~$\theta_0=0$.
Diffusions~\ref{e.mainSDE} and~\ref{e.one-arm} are identical when
$\kappa=6$, and \dref{e.mainSDE} (for~$4<\kappa<8$) is given
by \dref{e.one-arm} (for~$4<\kappa<\infty$) upon substituting
$\kappa\to(2\kappa)/(\kappa-4)$ and scaling time by $t\to
t(\kappa-4)/2$.  However, \dref{e.mainSDE} is a more singular
Bessel-type process when~$\kappa\leq 4$, requiring additional technical
analysis to deal with the times at which process is at~$0$
(see e.g.\ \lref{martingale}).  Furthermore, only the
large-time asymptotic decay rate of the hitting time distribution
(used in the proof of \tref{t.dimension}) is given in
\cite{LSW:one-arm}, and additional effort is required to obtain the
precise hitting time distribution provided in \tref{Bklaw}.

\subsection{Local martingales}

Recall the hypergeometric function defined by
$$ F(a,b;c;z) = \sum_{n=0}^\infty \frac{(a)_n (b)_n}{(c)_n n!} z^n \,,$$
where $a,b,c\in\C$ are parameters, $c\notin-\N$
(where $\N=\{0,1,2,\dots\}$), and $(\ell)_n$ denotes $\ell
(\ell+1)\cdots(\ell+n-1)$.  This definition holds for~$z\in\C$ when~$|z|<1$,
and it may be defined by analytic continuation elsewhere
(though it is then not always single-valued).  We define
for~$\lambda\in\C$
\begin{align*}
M^e_{\kappa,\lambda}(\theta)&=
  \textstyle F\!\left(1-\frac4\kappa+\sqrt{(1-\frac4\kappa)^2+\frac{8\lambda}{\kappa}},1-\frac4\kappa-\sqrt{(1-\frac4\kappa)^2+\frac{8\lambda}{\kappa}}; \frac32-\frac4\kappa; \sin^2\frac{\theta}{4}\right)\\
M^o_{\kappa,\lambda}(\theta)&=
  \textstyle F\!\left(1-\frac2\kappa+\sqrt{(\frac12-\frac2\kappa)^2+\frac{2\lambda}{\kappa}},1-\frac2\kappa-\sqrt{(\frac12-\frac2\kappa)^2+\frac{2\lambda}{\kappa}}; \frac32; \cos^2\frac{\theta}{2}\right) \cos\frac{\theta}{2}
\end{align*}
(where the formula for $M^e_{\kappa,\lambda}(\theta)$ makes sense
whenever $\kappa\neq \frac83,\frac85,\frac87,\dots$).  There is some
ambiguity in the choice of square root, but since
$F(b,a;c;z)=F(a,b;c;z)$, as long as the same choice of square root is
made for both occurences, there is no ambiguity in these definitions.

\begin{lemma} \label{l.localmart}
  For the diffusion \eqref{e.mainSDE} with~$\kappa>0$,
  let~$T$ be the first time at which
  $\theta_t\in2\pi\Z$, and let $\bar t=\min(t,T)$.  For any~$\lambda\in\C$,
  both $\exp[\lambda \bar t]M^e_{\kappa,\lambda}(\theta_{\bar t})$ and
  $\exp[\lambda \bar t]M^o_{\kappa,\lambda}(\theta_{\bar t})$ are
  local martingales parameterized by~$t$.
\end{lemma}

\begin{proof}
Given these formulas, in principle it is straightforward to verify
that the~$dt$ term of the \Ito\ expansion of
$$d \left[e^{\lambda t} M^{e|o}_{\kappa,\lambda}(\theta_t)\right]$$
is equal to zero (where $e|o$ is either $e$ or~$o$).
This term can be expressed as
$$
 e^{\lambda t}\left[\lambda\, M^{e|o}_{\kappa,\lambda}(\theta_t) + \left({M^{e|o}_{\kappa,\lambda}}\right)'(\theta_t) \,\frac{\kappa-4}{2} \,\cot(\theta_t/2) + \frac{\kappa}{2} \left({M^{e|o}_{\kappa,\lambda}}\right)''(\theta_t)\right]dt\,.
$$
Since Mathematica does not simplify this to zero, we show how to do this for $M^e_{\kappa,\lambda}$; the case for $M^o_{\kappa,\lambda}$ is similar.

We abbreviate $M^e_{\kappa,\lambda}$ with~$M$, let $F$ denote the hypergeometric function in the definition of
$M^e_{\kappa,\lambda}$,
let $a,b$, and $c$ denote the parameters of~$F$ in~$M$, and
change variables to $y=y(\theta)=\sin^2(\theta/4)=(1-\cos(\theta/2))/2$:
\begin{align*}
M(\theta) &= F(y(\theta))\\
M'(\theta) &= \frac14 F'(y(\theta))\sin(\theta/2)\\
M''(\theta) &= F''(y(\theta)) \frac1{16} \sin^2(\theta/2) + F'(y(\theta)) \frac{1}{8}\cos(\theta/2)\\ &= F''(y)\frac{y-y^2}{4} + F'(y)\frac{\frac12-y}{4}\,,
\end{align*}
so that
\begin{align*}
e^{\lambda t} \left[\lambda F(y) -\frac{\kappa-4}{4} F'(y) (y-{\textstyle\frac12}) + \frac\kappa8 F''(y) (y-y^2) + \frac\kappa8 F'(y)(\textstyle\frac12-y)\right] \\
 = e^{\lambda t} \left[\lambda F(y) + (1-3\kappa/8) F'(y) (y-{\textstyle\frac12}) + \frac\kappa8 F''(y) (y-y^2)\right] \\
 = e^{\lambda t} \sum_{n=0}^\infty y^n \frac{(a)_n (b)_n}{(c)_n n!}
\underbrace{\left[\begin{aligned}\lambda &+ (1-3\kappa/8)\left(n-\frac12\frac{(a+n)(b+n)}{c+n}\right)+{} \\&{}+ \frac\kappa8 \left(\frac{(a+n)(b+n)}{c+n} n -n(n-1)\right)\end{aligned}\right]}_{E_n/(c+n)}
\end{align*}
and we define $E_n$ to be $c+n$ times the expression in brackets.
(Note that indeed $c\notin-\N$.) We may write~$E_n$ as a polynomial
in~$n$:
\begin{align*}
E_n =
&[(1-3\kappa/8) (1-1/2) + (\kappa/8) (1-c+a+b)] \times n^2 +\\
&[\lambda + (1-3\kappa/8) (c-a/2- b/2) + (\kappa/8) (c+a\,b)] \times n +\\
&[c\lambda - (1-3\kappa/8) a\, b/2]\,.
\end{align*}
By our choices of $a$, $b$ and $c$, $E_n=0$ for each~$n$, which
proves the claim for $M^e_{\kappa,\lambda}$.
\ifcmp{\qed}
\end{proof}

\subsection{Expected first hitting of \texorpdfstring{$2\pi\Z$}{2\003\300 Z}}
\label{firsthit}

In this subsection we obtain asymptotics for the function
$$ L(\theta) := \E[\theta_T|\theta_0=\theta]\,,$$
where~$\theta_t$ is the diffusion~\eqref{e.mainSDE} and~$T$ is the
first time~$t\geq0$ at which $\theta_t\in 2\pi\Z$.  (Recall that
$T$ is finite a.s.\ when~$\kappa<8$.)  Whenever
$\theta\in2\pi\Z$, trivially $L(\theta)=\theta$.
\begin{lemma}
  For the diffusion~\eqref{e.mainSDE} with $8/3<\kappa<8$,
 $L(\theta_t)$ is a martingale.
\end{lemma}
\begin{proof}
  Since $L(\theta)$ is defined in terms of expected values,
  $L(\theta_t)$ is a local martingale whenever $\theta_t\notin2\pi\Z$,
  and the stopped process $L(\theta_{\min(t,T)})$ is a martingale.
  Since the diffusion behaves symmetrically around the points~$2\pi\Z$
  and the number of intervals of $\R\smallsetminus 2\pi\Z$ crossed
  before some deterministic time has exponentially decaying tails
  (which implies integrability), $L(\theta_t)$ is a martingale.
\ifcmp{\qed}
\end{proof}

Next, we express $L(\theta_t)$ in terms of the
$\lambda=0$ case of the local martingales $\exp[\lambda t]
M^e_{\kappa,\lambda}(\theta_t)$ and $\exp[\lambda t]
M^o_{\kappa,\lambda}(\theta_t)$.  Because $M^e_{\kappa,0}(\theta)=1$,
this local martingale is uninformative, but
$$
M^o_{\kappa,0}(\theta)=
  \textstyle F\left(\frac32-\frac4\kappa,\frac12; \frac32; \cos^2\frac{\theta}{2}\right)
  \cos\frac{\theta}{2}\,.
$$
We have $M^o_{\kappa,0}(0)=-M^o_{\kappa,0}(2\pi)=\sqrt{\pi}\Gamma(4/\kappa-1/2)/(2\Gamma(4/\kappa))$, and since $M^o_{\kappa,0}$ is bounded (when~$\kappa\neq8$), the stopped process $M^o_{\kappa,0}(\theta_{\min(t,T)})$ is also a martingale.
This determines~$L$, namely,
\begin{equation}
\label{Ltheta1}
  L(\theta) = \textstyle \pi -  \frac{2\,\sqrt\pi\,\Gamma(\frac4\kappa)}{\Gamma(\frac4\kappa-\frac12)} F\left(\frac32-\frac4\kappa,\frac12; \frac32; \cos^2\frac\theta2\right) \cos\frac\theta2\,,\qquad \theta\in[0,2\,\pi]
\end{equation}
and $L(\theta+2\,\pi)=L(\theta)+2\,\pi$. (It is also possible to
derive \eqref{Ltheta1} from the formula for $\Pr[\text{SLE trace
passes to left of $x+iy$}]$ \cite{MR1871700} after applying a
M\"obius transformation and suitable hypergeometric identities.)

\bigskip

We wish to understand the behavior of~$L$ near the points in~$2\pi\Z$,
and to this end we use the formula (see
\cite[pg.\ 108, Eq.~2.10.1]{\BE})
\begin{equation}
\label{e.Ftr}
\begin{aligned}
F(a,b;c;z) & = \frac{\Gamma(c)\Gamma(c-a-b)}{\Gamma(c-a)\Gamma(c-b)} F(a,b;a+b-c+1;1-z) + {} \\
& \qquad {}+ \frac{\Gamma(c)\Gamma(a+b-c)}{\Gamma(a)\Gamma(b)}(1-z)^{c-a-b} F(c-a,c-b;c-a-b+1;1-z)
\end{aligned}
\end{equation}
which is valid when $1-c$, $b-a$, and $c-b-a$ are not integers and $|\arg(1-z)|<\pi$.
 In our case the nonintegrality condition is satisfied whenever $8/\kappa\notin\N$.  For the range of~$\kappa$ that we are interested in, this rules out~$\kappa=4$, for which we already know $L(\theta)=\theta$, and therefore do not need asymptotics.  The endpoints of the range, $\kappa=8/3$ and $\kappa=8$ are also ruled out,
but for the remaining~$\kappa$'s we have
\begin{multline*}
\textstyle F(\frac32-\frac4\kappa,\frac12;\frac32;\cos^2\frac\theta2) =
 \displaystyle\frac{\Gamma(\frac32)\Gamma(\frac4\kappa-\frac12)}{\Gamma(\frac4\kappa)\Gamma(1)}
 \textstyle F(\frac32-\frac4\kappa,\frac12;\frac32-\frac4\kappa;\sin^2\frac\theta2) +{}\\
{}+ \frac{\Gamma(\frac32)\Gamma(\frac12-\frac4\kappa)}{\Gamma(\frac32-\frac4\kappa)\Gamma(\frac12)}
\textstyle\left|\sin\frac\theta2\right|^{\frac8\kappa-1}
 F(\frac4\kappa,1;\frac4\kappa+\frac12;\sin^2\frac\theta2)
\\
= \displaystyle\frac{\sqrt\pi\Gamma(\frac4\kappa-\frac12)}{2 \Gamma(\frac4\kappa)}
 \frac1{|\cos\frac\theta2|} 
+ \frac{\Gamma(\frac12-\frac4\kappa)}{2\Gamma(\frac32-\frac4\kappa)}
\textstyle\left|\sin\frac\theta2\right|^{\frac8\kappa-1}
 F(\frac4\kappa,1;\frac4\kappa+\frac12;\sin^2\frac\theta2)
\,,
\end{multline*}
and by~\eqref{Ltheta1}
$$
L(\theta) =
c_\kappa 
\textstyle\left(\sin\frac\theta2\right)^{\frac8\kappa-1}
 F(\frac4\kappa,1;\frac4\kappa+\frac12;\sin^2\frac\theta2) \cos\frac\theta2\,,
\qquad \theta\in(0,\pi)\,.
$$
for some constant~$c_\kappa>0$.
(One can use the Legendre duplication formula to show that
$c_\kappa= 2^{8/\kappa-1}{\Gamma(\frac4\kappa)^2}/{\Gamma(\frac8\kappa)}$, but we do not need this.)
Since $L(-\theta)=-L(\theta)$, we conclude that
$$
L(\theta)=A_0(\theta^2)\,|\theta|^{8/\kappa}/\theta\,,\qquad\theta\in(-\pi,\pi)\,,
$$
where $A_0$ is a analytic function (depending on~$\kappa$)
satisfying~$A_0(0)> 0$.
This implies
\begin{equation} \label{L(small)}
\theta^2 = A\bigl(|L(\theta)|^{2\kappa/(8-\kappa)}\bigr)
\end{equation}
near~$\theta=0$, for some analytic~$A$.

\subsection{Starting at \texorpdfstring{$\theta_0=0$}{\003\270=0}}

We will eventually need to start the diffusion at~$\theta_0=0$, but
\lref{l.localmart} only covers what happens up until the first
time that $\theta_t\in2\pi\Z$.  In this subsection we show
\begin{lemma} \label{martingale}
  For the diffusion \eqref{e.mainSDE} with $8/3<\kappa<8$,
  let~$T$ be the first time at which
  $\theta_t\in 2\pi\Z \smallsetminus 4\pi\Z$, and let $\bar
  t=\min(t,T)$.  For any $\lambda\in\C$, the process $\exp[\lambda \bar t]
  M^{e}_{\kappa,\lambda}(\theta_{\bar t})$ is a martingale.
\end{lemma}
\begin{proof}
Let $M$ abbreviate $M^{e}_{\kappa,\lambda}$, and let us assume without loss of
generality that $-2\pi < \theta_0 < 2\pi$.
Let us define $\omega_t=L(\theta_t)$, which is a martingale and may
be interpreted as a time-changed Brownian motion.
We wish to argue that $e^{\lambda {\bar t}} M(L^{-1}(\omega_{\bar t}))=
  e^{\lambda \bar t} M(\theta_{\bar t})$ is a local martingale.
(The definition of~$L$ implies that it is strictly monotone, and
hence $L^{-1}$ is well defined.) Note that by
\lref{l.localmart} it is a local martingale when $\theta_{\bar
t}\notin 2\,\pi\,\Z$. To extend this to a neighborhood of
$\theta_{\bar t}=0$, one could try to use It\^o's formula. To do
this, it would be necessary that $f:=M\circ L^{-1}$ be twice
differentiable. We have $M(\theta)=A_1(\theta^2)$ in
$(-2\,\pi,2\,\pi)$, where $A_1$ is analytic.
Consequently,~\eqref{L(small)} gives for $8/3<\kappa<8$ and for
$\omega$ in a neighborhood of~$0$,
\begin{equation}\label{fexpansion}
f(\omega)=
A_2\bigl(|\omega|^{2\kappa/(8-\kappa)}\bigr)\,,
\end{equation}
for some analytic~$A_2$.
Though this is not necessary for the proof, one can check that $A_2'(0)\ne 0$ and therefore $f''(0)$ is not finite
when~$\kappa<4$.

To circumvent the problem of $f=M\circ L^{-1}$ not being twice
differentiable, we use the It\^o-Tanaka Theorem (\cite[Theorem~1.5
on pg.~223]{\RY}). The exponent $\frac{2 \kappa}{8-\kappa}$
in~\eqref{fexpansion} ranges from $1$ to~$\infty$ as $\kappa$ ranges
from $8/3$ to~$8$. In particular, $f'(0)=0$ and $f'$ is continuous
near~$0$. Since $A_2$ is analytic, near~$0$ the function~$f$ may be
expressed as the difference of two convex functions, namely,
$f(\omega)=\bigl(f(\omega)-f(0)\bigr)\,1_{\omega\ge
0}+\bigl(f(\omega)-f(0)\bigr)\,1_{\omega\le 0}+f(0)$. Therefore, we
may apply the It\^o-Tanaka Theorem to conclude that $e^{\lambda \bar
t}f(\omega_{\bar t})= e^{\lambda \bar t} M(\theta_{\bar t})$ is a
local martingale also when $\theta_{\bar t}$ is near zero.

Now, the hypergeometric function~$F$ satisfies
\begin{equation} \label{F(1)}
F(a,b;c;1) = \frac{\Gamma(c)\Gamma(c-a-b)}{\Gamma(c-a)\Gamma(c-b)}
\end{equation}
provided $-c\notin\N$ and $\Re c>\Re(a+b)$ (see e.g.~\cite[pg.~104 eq.~46]{bateman-erdelyi:vol1}).
Therefore, $M(\pm 2\,\pi)$ is finite. Thus, $e^{\lambda \bar t} M(\theta_{\bar t})$ is bounded for bounded~$t$,
and we may conclude that it is a martingale.
\ifcmp{\qed}
\end{proof}

For future reference, we now calculate $M^{e}_{\kappa,\lambda}(\pm 2\,\pi)$.
Observe that the parameters $a,b,c$ of the hypergeometric function in the definition of
$M^e_{\kappa,\lambda}$ satisfy $2\,c-a-b=1$.
Consequently, the identity $\Gamma(z)\Gamma(1-z)=\pi/\sin(\pi z)$ and~\eqref{F(1)}
give
\begin{equation}\label{Mbd}
M^{e}_{\kappa,\lambda}(\pm 2\,\pi)=
\frac
{\sin\left(\frac\pi 2-\pi \sqrt{(1-\frac4\kappa)^2+\frac{8\lambda}{\kappa}}\right)}
{\sin\left(3\pi/2-4\pi/\kappa \right)}
=
  \frac{\cos\left(\pi\sqrt{(1-4/\kappa)^2+8\lambda/\kappa}\right)}{\cos(\pi(1-4/\kappa))}\,.
\end{equation}

\section{Proofs of main results}
We now restate and prove \tref{Bklaw}.
\begin{theorem} \label{t.mgf}
  Suppose the diffusion process \eqref{e.mainSDE} (with
  $8/3<\kappa<8$) is started at~$\theta_0=0$, and $T$ is the first
  time at which $\theta_t=\pm 2\pi$.  If $\Re\lambda\leq 0$, then
$$\E\bigl[e^{\lambda T}\bigm|\theta_0=0\bigr] = \frac{\cos(\pi(1-4/\kappa))}{\cos\left(\pi\sqrt{(1-4/\kappa)^2+8\lambda/\kappa}\right)}\,.$$
\end{theorem}
(This is equivalent to \tref{Bklaw} by
\pref{p.bklaw} and the remarks following the statement of
\tref{Bklaw}.)

\begin{proof}
Since
$M^e_{\kappa,\lambda}(\theta_T)=M^e_{\kappa,\lambda}(\pm2\,\pi)=
M^e_{\kappa,\lambda}(2\,\pi)$ a.s.\ and $\exp[\lambda {\bar t}]
M^{e}_{\kappa,\lambda}(\theta_{\bar t})$ is a martingale, the
optional sampling theorem gives
$$
M^e_{\kappa,\lambda}(2\,\pi) \,
\E\bigl[e^{\lambda T}\bigm|\theta_0=0\bigr] =
\E\bigl[e^{\lambda T}M^e_{\kappa,\lambda}(\theta_T)\bigm|\theta_0=0\bigr]
=
{M^e_{\kappa,\lambda}(0)}
=
1\,,
$$
and the proof is completed by appeal to~\eqref{Mbd}.
\ifcmp{\qed}
\end{proof}

\ifarxiv{\begin{proof}[Proof of \tref{t.dimension}]}
\ifcmp{
\newenvironment{pfof}[2]{\removelastskip\vspace{6pt}\noindent
 \textit{#1.}~#2}{\par\vspace{6pt}}
\begin{pfof}{Proof of \tref{t.dimension}}}
  Fix some~$\eps>0$ and
  let~$z\in\D$. Set $r_0:=1-|z|$, $r_1:=\dist(z,L_1^z)$, and suppose that $\eps<r_0$.
  We seek to estimate the probability that the open disk of radius~$\eps$ about~$z$ intersects the gasket; that is,
  the probability that $r_1<\eps$.
By the Koebe $1/4$ theorem, $r_1 \le \crad(D_1,z)\le 4\,r_1$.
Likewise, $r_0\le \crad(\D,z)\le 4\,r_0$.
Thus, $B_1^z=\log\crad(\D,z)-\log\crad(D_1,z)=\log(r_0/r_1)+O(1)$.
    Referring to the density function
  of $B_k^z$ \eqref{pdf}, we see that $$\Pr[r_1<\eps]\asymp \exp[-\alpha \log(r_0/\eps)] = (\eps/r_0)^\alpha\,,$$
  where
  $$\alpha = \frac{1/4-(1-4/\kappa)^2}{8/\kappa} = -\frac{3\kappa}{32} + 1 -\frac{2}{\kappa} = \frac{(8-\kappa) (3\kappa-8)}{32\kappa}\,.
  $$

For each $j=1,\dots,\lceil 1/\eps\rceil$, we may cover the annulus
$\bigl\{z:(j-1)\,\eps\le 1-|z|\le j\,\eps\bigr\}$ by $O(1/\eps)$
disks of radius~$\eps$. The total expected number of these disks
that intersect the gasket is at most
$$\sum_{j=1}^{\lceil 1/\eps\rceil} O(1/\eps) \times
  O(\eps/(j\eps))^\alpha = O(\eps^{\alpha-2})\,.$$
  (Here we made use of
  the fact that~$\alpha<1$.)  Thus on average $O(\eps^{\alpha-2})$
  disks of radius~$\eps$ suffice to cover the gasket.

  On the other hand, we may pack into~$\D$ at least $\Theta(1/\eps^2)$
  points so that every two of them are more than distance~$4\eps$
  apart, and each of them is at least distance~$1/2$ from the
  boundary.  For each such point~$z$ there is a $\Theta(\eps^\alpha)$
  chance that the disk or radius~$\eps$ centered at~$z$ is not
  surrounded by a loop, i.e., that that the gasket~$\Gamma$ contains a
  point~$z'$ that is within distance~$\eps$ of~$z$.  Since the points~$z$
  are sufficiently far apart, the points~$z'$ must be covered by
  distinct disks in any covering of~$\Gamma$ by disks of radius~$\eps$.
  Thus the expected number of disks of radius~$\eps$ required
  to cover the gasket~$\Gamma$ is at least $\Theta(\eps^{\alpha-2})$.
\ifcmp{\qed\end{pfof}}
\ifarxiv{\end{proof}}

\section{Open problems}

Kenyon and Wilson \cite{KW} also predicted the large-$k$ limiting
distribution of another quantity, the ``electrical thickness'' of
the loops~$L_k^z$ when $k\to\infty$.  The electrical thickness of a
loop compares the conformal radius of the loop to the conformal
radius of the image of the loop under the map $m(w)=1/(w-z)$, and
more precisely it is
$$\vartheta_z(L_k^z)=-\log\crad(L_k^z,z)-\log \crad(m(L_k^z),z)\,.$$
Kenyon and Wilson \cite{KW} predicted that the large-$k$ moment generating
function of~$\vartheta_z(L_k^z)$ is
\begin{equation} \label{th-mgf}
  \lim_{k\to\infty}\E[\exp(\lambda \vartheta_z(L_k^z))] = \frac{\sin(\pi(1-4/\kappa))}{\pi(1-4/\kappa)}\frac
         {\pi\sqrt{(1-4/\kappa)^2+ 8 \lambda /\kappa}}
         {\sin\left(\pi\sqrt{(1-4/\kappa)^2+ 8 \lambda /\kappa}\right)}\,,
\end{equation}
or equivalently that the limiting probability density function is
given by the density function of the exit time of a standard Brownian
\textit{excursion\/} started in the middle of the interval
$(-2\pi/\sqrt\kappa,2\pi/\sqrt\kappa)$, reweighted by a
factor of $\text{const}\times\exp[(\kappa-4)^2 x/(8\kappa)]$.
(This equivalence follows from
\cite[eq.~5.3.0.1]{borodin-salminen}.)
Recall that the density function of~$B_k^z$ is given by the density function
of the exit time of a standard Brownian \textit{motion\/} started in
the middle of the interval $(-2\pi/\sqrt\kappa,2\pi/\sqrt\kappa)$, also
reweighted by a factor of $\text{const}\times\exp[(\kappa-4)^2
x/(8\kappa)]$.  These forms are highly suggestive, but currently we do
not know how to calculate the electrical thickness using \CLEk/, nor do
we have a conceptual explanation for why these distributions take
these forms.

\pdfbookmark[1]{References}{bib}
\bibliographystyle{hmralpha}
\bibliography{cle-radii}

\end{document}